\documentclass{aupl}

\usepackage{
amsfonts,
latexsym,
amssymb,
enumerate,
verbatim,
mathrsfs,
graphicx,
centernot,
}

\newcommand{\labbel}{\label}

\newtheorem{theorem}{Theorem}[section]

\newtheorem{proposition}[theorem]{Proposition} 
 
\newtheorem{corollary}[theorem]{Corollary}

\newtheorem*{claim*}{Claim}

\newtheorem*{theorem*}{Theorem}
\newtheorem*{proposition*}{Proposition}
\newtheorem*{corollary*}{Corollary}
\newtheorem*{lemma*}{Lemma}
\newtheorem*{scholion*}{Scholion}

\theoremstyle{definition}
\newtheorem{definition}[theorem]{Definition}

\newtheorem{problem}[theorem]{Problem}

\theoremstyle{remark}
\newtheorem{remark}[theorem]{Remark}

\newtheorem*{remark*}{Remark}
\newtheorem*{remarks*}{Remarks}
 
\newtheorem{example}[theorem]{Example}

\newtheorem{observation}[theorem]{Observation}

\newtheorem*{observation*}{Observation}

 \allowdisplaybreaks[1]

\numberwithin{equation}{section}

\begin{document}

\title
{Varieties defined by 
linear equations  have the  amalgamation property}

\author{Paolo Lipparini} 
\address{Dipartimento di Matematica\\Variet\`a  della  Ricerca
 Scientifica\\Universit\`a di Roma ``Tor Vergata'' 
\\I-00133 ROME ITALY}

\email{lipparin@axp.mat.uniroma2.it}

\subjclass{03C05;  03C52; 08B05;  08B25}

\keywords{Linear equation; equilinear equation; variety;
amalgamation property; strong amalgamation property;
joint embedding property; Maltsev condition; Fra\"\i ss\'e limit}

\thanks{\today\\
Work performed under the auspices of G.N.S.A.G.A. Work 
partially supported by PRIN 2012 ``Logica, Modelli e Insiemi''.
The author acknowledges the MIUR Department Project awarded to the
Department of Mathematics, University of Rome Tor Vergata, CUP
E83C18000100006.}

\begin{abstract}
A variety is a class of algebraic structures axiomatized
by a set of equations.
An equation is \emph{linear}
if there is at most one occurrence of an operation  symbol
on each side.
We show that a variety axiomatized by linear equations
 has  the strong
amalgamation property. 

Suppose further that  the language has no constant symbol 
and, for each equation,
either  one side is operation-free, or
exactly the same variables appear
on both sides. Then also
 the joint embedding property holds.

Examples include most varieties defining classical Maltsev conditions. 
In a few special cases, the above properties 
are preserved when  further unary operations
appear in the equations.
 \end{abstract}

\maketitle

\section{Introduction} \labbel{intro} 

An \emph{equation} is an atomic formula in a language
with equality as the only relational symbol. 
We call an equation $\sigma$  \emph{linear} if 
$\sigma$ has at most one occurrence of  operation symbols on each side.
Constants are not considered as operations.
If in addition the set of variables appearing on the right is
equal to the set of variables appearing on the left,
we say that $\sigma$ is \emph{equilinear}. 
We include in the class of equilinear equations
also the equations in which only one operation occurs,
that is, one side of the equation is a variable or a constant.

Equilinear equations are important
in universal algebra
 since the great majority of classical Maltsev
 conditions are obtained  by interpreting
a variety   axiomatizable by
equilinear equations.
Recall that a \emph{variety}
is a class axiomatized by a set of equations
in a language with equality as the only relation.
 Every nontrivial
idempotent Maltsev condition
implies a nontrivial equilinear idempotent Maltsev condition.
See \cite{Ta},  \cite[Chapter 9]{HMK}  and \cite[Section 2.4]{KK}. 
In passing, it is worth mentioning that there exists
the weakest idempotent Maltsev condition \cite{Sig,Ols},
though this surprising recent result will play no role here. 

As an argument
more familiar to model theorists,
the amalgamation property \cite{F} is an important notion
which has found many successful applications 
 also in algebra, algebraic logic and  category theory
\cite{Ev,GM,H,J,KMPT,MMT}.
From the perspective of universal algebra,
the importance of the property stems to the fact that
a variety $\mathcal V$ has the amalgamation property
if and only if the category of algebras in $\mathcal V$
with embeddings has push outs \cite{KMPT}.
Here \emph{algebra} is a shorthand, for \emph{algebraic structure},
namely, a model in a language without relation symbols.  

In Section \ref{gen} we show that a variety axiomatizable by linear equations,
possibly containing constants,
 has  the strong
amalgamation property.
The case of a  variety axiomatizable by equilinear equations
in a language without constants
is much simpler and is treated in Section \ref{ap}. 
In this case we 
also get the joint embedding property.
At large, the above results provide a striking partial 
affirmative answer
to an ancient problem by B. J{\'o}nsson, asking
for \emph{general  results  that  assert  that  if  an  elementary  
class  is  characterized 
by  axioms of such and  such a form,  then  this class has  the  amalgamation 
property} \cite{J}.

The above results apply to virtually all the 
varieties used to define classical Maltsev conditions,
such as the variety 
with a Maltsev operation, the variety 
with $n$ J{\'o}nsson  operations, the variety  
with $n$ Day operations, 
 etc. 
Since we show that also the class of finite algebras
in such varieties  has
both the amalgamation property and the joint embedding
property, a classical result by Fra\"\i ss\'e \cite{F}  applies, to the
effect that such classes have 
a Fra\"\i ss\'e limit, 
an ultrahomogeneous universal countable
algebra. See \cite[Theorem 7.1.2]{H}.
Such ``generic'' algebras probably deserve
an accurate study, but we shall not deal with this here. 

We conclude the note with a few remarks in Section \ref{add}  on adding
unary operations.
If $\mathcal V$ is a variety in the language $\mathcal L$, let $\mathcal V_h$
be a variety in the  language $\mathcal L \cup \{  h\} $,
where $h$ is a new unary function symbol.
 It is almost trivial that if $\mathcal V$ 
is a variety with the amalgamation property, then 
$\mathcal V_h$ has the amalgamation property, if  $h$ satisfies
the  equations asserting that 
$h$ is a homomorphism. 
A counterexample is provided showing that
the amalgamation property is not necessarily 
preserved if $h$ is assumed to respect only a subset
of the operations of $\mathcal V$. 
However, the amalgamation property is preserved in $\mathcal V_h$ 
when $h$ has an inverse and $\mathcal V$ is axiomatized by  
 linear equations, no matter which subset of
operations is respected by $h$.

\section{Preliminaries} \labbel{prel}

We give the definition of the 
strong amalgamation property in  the case of algebraic structures.
The general model-theoretical definition presents no significant
difference.
For simplicity, we assume to deal with a class  $\mathcal K$ of algebras
such that  $\mathcal K$ is
closed under taking isomorphism.
This is no big loss of generality; in fact, we
shall almost invariably deal with varieties.  

\begin{definition} \labbel{apdef}   
A class  $\mathcal K$ closed under  isomorphism
is said to have the \emph{strong amalgamation property},
SAP, for short,
if the following holds.
For every three algebras $\mathbf A$, $\mathbf  B$ and $\mathbf  C$
in $\mathcal K$ 
such that 
\begin{equation}\labbel{apeq}
\text{$\mathbf  C \subseteq \mathbf A$, \quad 
 $\mathbf  C \subseteq \mathbf B$ \quad 
and \quad  $A \cap B = C$,} 
   \end{equation}
there is an 
\emph{amalgamating} algebra,
or an \emph{amalgam}, for short,   $\mathbf  D$ in $\mathcal K$  such that    
$\mathbf  A \subseteq \mathbf D$ and 
 $\mathbf  B \subseteq \mathbf D$.
Here, say, $\mathbf  C \subseteq \mathbf A$
means that $C \subseteq A$ as sets, and  $\mathbf  C $ is
a substructure of $  \mathbf A$. 

We shall  usually prove  the strong amalgamation
property; however, in some comments we shall deal
with the plain (=not strong) amalgamation property.
A class  $\mathcal K$ closed under  isomorphism
is said to have the \emph{amalgamation property},
AP, for short, if, under the assumptions in 
\eqref{apeq}, there 
 are an algebra $\mathbf  D \in \mathcal K$  
and embeddings
$\iota : \mathbf A \to \mathbf  D$ 
and
$ \kappa  : \mathbf B \to \mathbf  D$
which agree on $C$.

The point is that  if the strong amalgamation property
holds, then
$\mathbf  D$ as above can be obtained 
in such a way that two elements $a \in A \setminus B$ 
and $b \in B \setminus A$
are never identified. On the other hand, it might happen that 
 $ \iota(a)= \kappa (b)$ for such elements,
if only AP is assumed. This is the reason why
in the definition of AP
we need to deal with embeddings, not simply
inclusions.

The class $\mathcal K$ has the
\emph{joint embedding property},
JEP, for short,
if, 
for every  $\mathbf A, \mathbf  B \in \mathcal K$,
there are an algebra $\mathbf  D \in \mathcal K$  
and embeddings
$\iota : \mathbf A \to \mathbf  D$ 
and
$ \kappa  : \mathbf B \to \mathbf  D$. 
 \end{definition}

We do not allow algebras with empty domain;
hence, in general,  JEP is not a consequence of AP.

In the present note we do not consider constants as operations, hence, in general, 
a linear equation is allowed to contain constants.
This  general case shall be treated in Section \ref{gen}.
We first deal with the simpler case when constants are not allowed.

\begin{definition} \labbel{lindef}    
A \emph{linear equation}
in a language without constants
 is an equation which
has one of the following forms. 
\begin{align} \labbel{line}  
f(x_{1},x_{2}, \dots, x_{n} ) &= y_1,  \text{ or } 
\\
\labbel{d}  
f(x_{1},x_{2}, \dots, x_{n} ) &= g(y_{1},y_{2}, \dots, y_{m}    ),
 \end{align}
where $f$ and $g$
are not necessarily distinct operations and 
$x_{1},x_{2}, \dots, y_{1},y_{2}, \dots$
are  not necessarily  distinct variables.
A more general definition of a linear equation
for languages with constants shall be given in Definition \ref{lingen}.

Notice that the terminology is not uniform in the 
literature. The above
terminology is  the most common
when dealing with Maltsev conditions. 

We say that an equation $\sigma$ is
\emph{equilinear} if either $\sigma$  has  the form \eqref{line},
or $\sigma$  has the form \eqref{d} and the set of variables appearing on the
 left-hand side
is equal to the set of variables appearing on the right-hand side, that is, 
$ \{ x_{1},x_{2}, \dots, x_{n} \}  = \{ y_{1},y_{2}, \dots, y_{m}  \} $,
with no consideration of multiplicities.
\end{definition}

Notice that an equilinear equation may logically 
imply a linear but not equilinear equation.
For example, the equilinear equation $x=f(x,y,y)$
(relevant to the first discovered Maltsev condition, see the comments
shortly after Corollary \ref{v})
logically implies the equation
$f(x,y,y)=f(x,z,z)$,
which is not equilinear.
 Hence we shall always need to refer to varieties $\mathcal V$  
\emph{axiomatized}
 by equations of some particular form, rather
than to varieties $\mathcal V$ such that 
 \emph{all the equations valid in $\mathcal V$}
 have some particular form.

\section{The strong amalgamation property for equilinear varieties} \labbel{ap}

 We first deal with languages without constants.
This is the case of most interest; in this case  
proofs are much simpler  and we get slightly stronger results.

\begin{theorem} \labbel{thm}
Suppose that  $\mathcal V$ 
is a variety 
axiomatizable by equilinear equations and  
the language  of $\mathcal V$ has no  constant.
Then $\mathcal V$ 
 has both the strong amalgamation property and 
 the joint embedding property. 
\end{theorem}

\begin{proof}
Suppose that
$\mathcal V$ is axiomatized by 
a set $\Sigma$  of equilinear equations.

First, notice that 
there might be further equations of the form \eqref{line} 
which follow from the given set of equations $\Sigma$.
For example, if $f(x,x,y,z)=g(x,x,y,z)$
and $g(x,x,y,y)=y$ are among the equations defining
$\mathcal V$, then  $f(x,x,y,y)=y$ holds in $\mathcal V$, as well,
even if this equation is not explicitly listed in $\Sigma$.
However, we can suppose without loss of generality to have added
to $\Sigma$  all the equations of the form  \eqref{line}
which follow logically from the given equations.
Henceforth we suppose to have extended  the  set $\Sigma$ in this way.

Furthermore, we can suppose that, for every equation in $\Sigma$
of the form \eqref{line},  $y_1 \in \{ x_{1},x_{2}, \dots  \}$,
since otherwise $\mathcal V$ is a trivial variety 
with only $1$-element algebras, and in this
case the conclusion is obvious.  

Given a triple $\mathbf A$, $\mathbf  B$, $\mathbf  C$
of algebras in $\mathcal V$ as in \eqref{apeq},
choose an arbitrary element  
$ \bar{d} \in A \cup B $.
For each operation $f $ in the language of $\mathcal V$,
 define $f$ on $A \cup B$ by
\begin{equation}\labbel{alf}    
 \begin{aligned} 
 f(d_1, d_2, \dots ) &= f _{\mathbf A} (d_1, d_2, \dots ) 
&&\text{if  $d_1, d_2, \dots \in A$,} 
\\
 f(d_1, d_2, \dots ) &= f _{\mathbf B} (d_1, d_2, \dots )
&&\text{if  $d_1, d_2, \dots \in B$,}
\\
  f(d_1, d_2, \dots ) &= d_j
&&
\begin{aligned}   
&\text{if this is forced by some equation}
\\[-2 pt]
 &\text{in $\Sigma$ of the form \eqref{line}, and} 
  \end{aligned}
\\
 f(d_1, d_2, \dots ) &= \bar{d}  
&&\text{if  none of the above clauses apply,} 
   \end{aligned}
\end{equation}   
for $d_1, d_2, \dots \in A\cup B$ and  where, in more detail, 
we apply the third clause if 
there is an equation in $\Sigma$ 
of the form $f(x_{1},x_{2}, \dots  ) = x_j$ 
 and the expression $ f(d_1, d_2, \dots )$
is obtained by substituting equal elements for equal variables,
 that is,
$d_h=d_k$, whenever  $x_h=x_k$. 
Notice that, as we mentioned, 
 $x_j \in \{ x_{1},x_{2}, \dots  \}$,
since otherwise we are in a trivial variety.
In particular, also $d_j \in \{ d_{1},d_{2}, \dots  \}$,
so that the definition makes sense.
Moreover, the definition is well-posed,
since if we have both 
$f(x_{1},x_{2}, \dots, x_{n} ) = y_1$
and
$f(x_{1},x_{2}, \dots, x_{n} ) = y_2$ in $\Sigma$,
for $y_1$ and  $y_2$ distinct variables, then we are
again  in a trivial
variety.    
 
With the above definition, $D=A \cup B$
becomes an algebra $\mathbf  D$ for the language of $\mathcal V$.
Notice that the first three clauses possibly overlap,
but they agree in any overlapping case, since $\mathbf  C$ is
a subalgebra of  
both  $\mathbf A$ and $\mathbf  B$; moreover,
the equations in $\Sigma$ of the form 
\eqref{line} hold in the three algebras to
be amalgamated.
By the first two
clauses in \eqref{alf}, $\mathbf  D$ extends both $\mathbf A$ and $\mathbf  B$.

It remains to show that the defining equations \eqref{line} and \eqref{d} 
are valid in $\mathbf  D$.
Since, by assumption, the equations are valid both in $\mathbf A$
 and in $\mathbf  B$,
 the conclusion is obvious if some
identity involves only elements from $A$ or only elements from $B$.
By the third clause 
in \eqref{alf}, $\mathbf  D$ satisfies all the equations in
$\Sigma$ of the form  \eqref{line}.
Hence it remains to consider an equation 
of the form  
\begin{equation}\labbel{id} 
    f(x_{1},x_{2}, \dots  ) = g(y_{1},y_{2}, \dots  )
  \end{equation}
from $\Sigma$, that is, we have to prove that 
$f(d_{1},d_{2}, \dots  ) = g(e_{1},e_{2}, \dots  )$,
when $d_{1},d_{2}, \dots, e_{1},e_{2},\dots \in A \cup B$
and equal elements are substituted for equal variables
in \eqref{id}. If the value of,  say, $f(d_{1},d_{2}, \dots  )$
is forced by some equation in $\Sigma$ of the form \eqref{line}, then the value of
$g(e_{1},e_{2}, \dots  )$ 
 is forced, too, and the resulting value is the same.
The reason for this fact is that  we have completed \eqref{line} with all the equations
logically following from all the equations defining the variety.

As  side remark, notice that the possibility of using the third clause
in combination with some equation of the form \eqref{id} 
does not necessarily trivialize the equation \eqref{id}.
For example, suppose that \eqref{id} is
\begin{equation}\labbel{pul}     
f(x,x,y,y,z) = g(x,x,y,y,z) 
  \end{equation}
and we need to check the instance 
$f(d,d,d,d,e) = g(d,d,d,d,e) $ of \eqref{pul}. 
It might happen that 
$f(x,x,x,x,y)=y$
is an equation in $\Sigma$,
so that  $f(d,d,d,d,e) \allowbreak =e$
is ``forced''.
Then
\eqref{pul} and  
$f(x,x,x,x,y)=y$ imply
that $g(x,x,x,x,y)=y$
is an equation in the  ``extended''   $\Sigma$,
hence  $g(d,d,d,d,e)=e$
is forced, as well.
However it might happen that, say, 
$f(d,d,d',d',e)$
is not forced,
hence neither  
$g(d,d,d',d',e)$
is forced.
If this is the case,
the definition in \eqref{alf}
provides   
$f(d,d,d',d',e) = \bar{d} = g(d,d,d',d',e)$,
a fact we are now going to exploit in general.

So far, we have seen that all the equations 
in $\Sigma$ of the form 
\eqref{line} are valid in $\mathbf  D$ 
and that, given an instance of \eqref{id},
the value of   the left-hand side is forced by the third clause
in \eqref{alf} 
if and only if the right-hand side is forced;
in this case, the two forced values are the same.
Hence we need  
to deal with the case when no side is ``forced''.
Since, by assumption, we are dealing with
an equilinear equation, 
we have  
$ \{ d_{1},d_{2}, \dots   \}  = \{  e_{1},e_{2}, \dots   \} $
in any instance of \eqref{id},
hence  $ \{ d_{1},d_{2}, \dots   \} \subseteq A $
if and only if  $\{  e_{1},e_{2}, \dots   \} \subseteq A$.
If this is the case,
$f(d_{1},d_{2}, \dots  ) $ and $  g(e_{1},e_{2}, \dots  )$
are given by the first clause in \eqref{alf}, and they are equal
since \eqref{id} is valid in $\mathbf A$.
A similar argument applies when
  $ \{ d_{1},d_{2}, \dots   \}  = \{  e_{1},e_{2}, \dots   \}  \subseteq B$.
In all the other cases
$f(d_{1},d_{2}, \dots  ) $ and $   g(e_{1},e_{2}, \dots  )$
are both given by the fourth clause in
\eqref{alf}, hence  
$f(d_{1},d_{2}, \dots  )= \bar{d}  
=  g(e_{1},e_{2}, \dots  )$.
We have checked that all the instances of \eqref{id}
hold in $\mathbf  D$, hence $\mathbf  D$
is an amalgamating algebra in $\mathcal V$.
The amalgamation property follows.

To prove the joint embedding property
assume, without loss of generality,
that $A \cap B = \emptyset $ and 
repeat the above arguments
considering $C= \emptyset $. 
 \end{proof} 

The assumption that at most one operation symbol occurs on each
side is obviously necessary in Theorem \ref{thm}.
The variety of semigroups is axiomatized by a single equation
  such that the same variables appear on both sides; however,
the variety of semigroups has not AP, a result credited to 
Kimura \cite{Ki} in \cite{J,KMPT}.

As we mentioned, we do not allow algebras  with 
empty domain.
Were we allowing empty
domains,  the \emph{strong} amalgamation property
would fail in the case of a trivial variety, 
taking $ C= \emptyset $
and $\mathbf A$ and $\mathbf  B$ $1$-element
algebras.
Obviously, the amalgamation property holds
even  in this extraordinary situation.

The algebra $\mathbf  D$ 
in the proof of Theorem \ref{thm}
has been constructed on the union of $A$ and $B$.
As we have discussed at length in \cite{apu},
this fact has interesting consequences, 
frequently with easy
proofs.  For example, 
items (2) and (4) in the following corollary
follow immediately from the above observation.

\begin{corollary} \labbel{v}
If $\mathcal V$ is a variety axiomatizable by a set
of linear equations,
then the following statements hold. 
  \begin{enumerate}   
 \item 
If $\mathcal V$ is not trivial, then, for every nonzero $n \in \mathbb N$,
$\mathcal V$ has an algebra with exactly $n$ elements.
\item
The class of finite algebras in $\mathcal V$ 
has the strong amalgamation property.
\item
If $\mathcal V$ is axiomatized by a set 
of equilinear equations and 
the language of $\mathcal V$ has no constant, then the
 class of finite algebras in $\mathcal V$ 
has the the joint
embedding property.
\item
Let $\Gamma$ be a set of (possibly infinitary)
universal-existential sentences for the language of $\mathcal V$ and
such that in each sentence 
at most one variable is bounded by the universal quantifier.

Then the class consisting of the algebras in $\mathcal V$ satisfying
all the   sentences in $\Gamma$
has the strong amalgamation property. 
  \end{enumerate} 
 \end{corollary} 

\begin{proof} 
We first give the proof
for varieties satisfying the assumptions in Theorem \ref{thm}.
The proof for the general case shall be postponed
to Section \ref{gen}. 
  
If $\mathcal V$ is not trivial, the proof of Theorem \ref{thm}
shows that an algebra $\mathbf  D$ witnessing JEP
can be obtained over the disjoint union of $A$ 
and $B$. Then (1) follows by induction on $n$,
since every variety has an algebra with one element. 
 The  statement in (3), too, follows from the above observation.
As we mentioned, all the rest follows from the fact that
the proof of \ref{thm} provides strong amalgamation
``into union''. 
\end{proof}  

Starting from A. Maltsev breakthrough 
characterization of congruence permutable varieties \cite{Ma},
Maltsev conditions have since played a predominant role in universal algebra.
Informally, a Maltsev condition asserts the existence of 
a finite set of terms satisfying a given finite set of equations.
More precisely, a Maltsev condition is determined by 
some finitely presented variety $\mathcal W$;
some variety $\mathcal V$  satisfies the Maltsev condition
associated to $\mathcal W$ if $\mathcal W$
is interpretable in   $\mathcal V$.
We shall not need the exact details here,
this discussion serves only as a motivation;
henceforth, we shall only briefly hint to the simplest example, the one discovered
by Maltsev himself. Here the name of Maltsev is attributed
both to the general notion and to the specific example.

The \emph{variety $\mathcal W_M$ with a Maltsev operation}
is the variety in the language with a ternary operation $f$ such that 
$\mathcal W_M$  is
axiomatized by the equations
\begin{equation*}\labbel{mal}
x=f(x,y,y), \qquad\qquad f(x,x,y)=y.
  \end{equation*}       
A variety $\mathcal V$
satisfies  the \emph{Maltsev condition} $\mathfrak M$ 
associated to $\mathcal W_M$
if $\mathcal V$ interprets  $\mathcal W_M$, namely,
 if   $\mathcal V$ has a ternary \emph{term} $t$ 
satisfying the equations   
$x=t(x,y,y) $, $  t(x,x,y)=y$. 

It is customary to deal with Maltsev conditions 
rather informally, expressing them, as above, by means of the 
existence of appropriate terms. However, here we need to
deal with varieties. Quite informally, one simply should 
replace the word ``term'' (sometimes, ``polynomial'' or
``expression'' 
in the classical literature) by ``operation''  
in, say, Maltsev \cite{Ma},  Pixley \cite{Pi},
Day \cite{Da}, J{\'o}nsson \cite{Jo},
Hagemann and Mitschke  \cite{HM}.

\begin{corollary} \labbel{cor}
The conclusions of Theorem \ref{thm}
and Corollary \ref{cor} 
 hold for the following varieties.
  \begin{enumerate}    
\item  
The  variety  with a Maltsev operation \cite{Ma}.
\item 
 The  variety  with a Pixley operation \cite{Pi}; and, for every n:
\item
The variety with an $n$-ary near unanimity operation.
\item
The  variety  with $n$ Day  operations \cite{Da}.
\item
The  variety  with $n$ J{\'o}nsson   operations \cite{Jo}.
\item
The  variety  with $n$ Hagemann-Mitschke  operations \cite{HM}.
 \end{enumerate} 
 \end{corollary} 

Corollary \ref{cor}
is just  illustrative.
Other examples of Maltsev conditions
such that Theorem \ref{thm} applies to their defining
varieties  can be found,
among many others and with partial overlaps, in  
Taylor \cite[Corollary 5.3]{Ta}, 
Gumm \cite[Theorem 7.4(iv)]{Gu},
Tschantz \cite[Lemmata 3 and 4]{Ts},  
Hobby, McKenzie \cite[Lemmata 9.4(3), 9.5(3),
Theorems 9.8(4), 9.11(4) and 9.15(3)]{HMK},
Siggers \cite[Theorem 1.1, Section 3]{Sig}, 
Kearnes, Kiss \cite[Theorems 3.21(3), 4.7, 5.23(3),   
5.28(3), 8.13(3) and 8.14(3)]{KK}, 
Kearnes,  Markovi\'{c}, McKenzie \cite[Theorem 2.2, Corollaries 3.1, 3.2]{KMM}, 
Ol\v{s}\'{a}k \cite[Theorems 3.3, 6.1, Definition 5.1]{Ols},
Kazda,  Kozik,   McKenzie, M. Moore \cite[Section 1]{KKMM}, 
Kazda, Valeriote \cite[Sections 3.2 and 3.3]{KV}, 
Lipparini \cite[Definitions 2.1, 2.7, 6.1, 7.1, 7.6, 8.1, Remarks 6.4, 8.16, 8.19, 10.11(c)]{daysh}.

Notice that in all the examples dealing with Maltsev 
conditions the defining variety is finitely presented, 
namely,  a variety axiomatized by a finite set of equations in a finite language.
On the other hand, no finiteness assumption is necessary in 
Theorem \ref{thm}.

Turning to another consequence
of Theorem \ref{thm}, we briefly recall
an important model theoretical construction. 
R. Fra\"\i ss\'e \cite{F} developed 
a by now classical method for obtaining 
some limit ``random'' or ``generic'' homogeneous and
universal structure,  starting from a countable set of finite
or countable structures 
with AP and JEP. This \emph{Fra\"\i ss\'e limit}
is the ``most general'' homogeneous countable structure
into which all the starting structures can be embedded.

The classical example is  the ordered set  of the rationals, which 
is, in a sense, the most general
countable linearly ordered set, and is the Fra\"\i ss\'e 
limit of the class of finite linear orders. 
As another example, the countable random graph
 is the limit of the class of finite graphs.
The random graph can be obtained in a number of different
ways; quite astonishingly, 
Fra\"\i ss\'e could have explicitly introduced it
(apparently, he did not) about a decade  
before it  became a subject of deep studies.
See \cite[Chapter 7]{H} for more details and precise definitions.

\begin{corollary} \labbel{frai}
Suppose that $\mathcal V$ is a variety in a finite language
without constants
and $\mathcal V$ is axiomatizable by equilinear equations. 
Then the class of finite algebras in $\mathcal V$ 
has a Fra\"\i ss\'e limit.

In particular, for every variety listed in Corollary \ref{cor}
and in the comment below it,
 finite algebras have a Fra\"\i ss\'e limit.
 \end{corollary} 

\begin{proof}
Immediate from
Corollary \ref{v}(2)(3) and 
 Fra\"\i ss\'e's Theorem; see, e.~g.,  \cite[Theorem 7.1.2]{H}.
\end{proof}

\begin{problem} \labbel{frpr}
It is probably interesting to study
the Fra\"\i ss\'e limits of the varieties listed in
Corollary \ref{cor}, as well as of the varieties referred to shorty after. 
 \end{problem}

\begin{remark} \labbel{wfin}    
The assumption that the language of $\mathcal V$ is finite
 in Corollary \ref{frai}
can be weakened to the assumption that $\mathcal V$
contains only a countable number of finite algebras modulo isomorphism. 

As an example to which this more general version of Corollary \ref{frai}
applies, consider the  variety $\mathcal V$ in a language
with operations $( f_n ) _{ n \in \mathbb N^+ } $,
each $f_n$ $n$-ary.  
The equations axiomatizing $\mathcal V$    are
\begin{equation} \labbel{*}  
f_{n+1}(x_1, x_2, \dots, x_{n+1}) =
f_{n}(y_1, y_2, \dots, y_{n}),
  \end{equation}    
for every $n >0$,  where at least one variable is repeated in the sequence
$(x_1, x_2, \dots, x_{n+1})$ and the sequence
$(y_1, y_2, \dots, y_{n})$ is obtained from 
$(x_1, x_2, \allowbreak \dots, \allowbreak  x_{n+1})$
by deleting the first occurrence of a repeated variable.

The equations \eqref{*}
are equilinear and, for every $n$,
there is only a finite number of algebras in $\mathcal V$
of cardinality $n$ 
up to isomorphism.
Indeed, in an algebra of cardinality $n$, all the operations of arity 
$>n$ are defined in function of     
$f_n$. 

Of course, the above remarks apply
to any subvariety of $\mathcal V$
axiomatized by equilinear equations. 
\end{remark}

\begin{remark} \labbel{varr}
(a) There are  possible  variations on the construction of $\mathbf  D$ 
in the proof of Theorem \ref{thm}.
Such variations might be of some interest; in fact, we shall see 
an application in the proof of Proposition \ref{un}  below.

(b) First of all, as we mentioned, there are interesting consequences
of the fact that we can choose $\mathbf  D$ to be over 
$D= A \cup B$. However, for nontrivial
varieties, the proof of Theorem \ref{thm}
works for every    $D \supseteq  A \cup B$, fixing 
any $ \bar{d}  \in D$.

(c) Moreover, there is no need to consider a fixed element 
$ \bar{d}$.  We get more varied amalgamating algebras 
if we fix a function
$\eta: D^{fin} \to D$, where
$D^{fin}$ is the set of all finite subsets of $D$ and 
$D \supseteq  A \cup B$.
Then replace the identity in the fourth clause in \eqref{alf}
by
$f(d_1, d_2, \dots, d_n)= \eta(\{ d_1, d_2, \dots, d_n \})$.
The proof of Theorem \ref{thm}
works with this definition, too, since 
when we 
prove that an identity of the form  \eqref{id} is valid in $\mathcal V$
and we 
apply the fourth clause in \eqref{alf}   
to 
$f(d_{1},d_{2}, \dots  ) $ and
$    g(e_{1},e_{2}, \dots  )$,
we 
have
  $ \{ d_{1},d_{2}, \dots   \}  = \{  e_{1},e_{2}, \dots   \} $.
Thus, under the present
more general definition, 
we have $f(d_{1},d_{2}, \dots  )=  \eta (\{ d_1, d_2, \dots \})=
\eta (\{ e_1, e_2, \dots \})
=  g(e_{1},e_{2}, \dots  )$.

(d) There are many ways to choose some function $ \eta$
as above. Of course,  we can pick a fixed element
$ \bar{d}  \in D$ and let $ \eta$ be the constant function
with value $ \bar{d}$. In this case, if $D= A \cup B$,  we get back the structure 
constructed in the proof of \ref{thm}.   
As another possibility, linearly order $D$ 
and let 
$ \eta(\{ d_1, d_2, \dots, d_n \}) = \max\{ d_1, d_2, \dots, d_n \}$. 

(e) Considering the possibility
$D \supset A \cup B$
provides also a slightly different proof for
item (1) in Corollary \ref{v} (so far, limited to 
the equilinear case without constants).
Take $\mathbf A=\mathbf  B=\mathbf  C$  
 $1$-element algebras and apply the 
proof of Theorem \ref{thm}, as modified above,
 by considering any set $D$ of cardinality $n$
and containing $C$.   
 \end{remark}

\section{The general case} \labbel{gen} 

We now treat the general case
when the language of $\mathcal V$ might contain constants.
We also deal with linear not necessarily
equilinear equations.
We still get the strong amalgamation property,
but the joint embedding property fails, in general,
as shown by the following examples.
See, however, Corollary \ref{corthmc}
for some special cases in which JEP holds
even in the presence of constants.

\begin{remark} \labbel{con}
(a) The assumption that the language 
of $\mathcal V$ has no constant is necessary
in Theorem \ref{thm}
in order to get the joint embedding property. 
For example, let the language of $\mathcal V$
have two constants $c_1$ and $c_2$.
If  $c_1$ and $c_2$ are interpreted in such a way
that $c_1 = c_2$ in $\mathbf A$ and $c_1 \neq c_2$
in $\mathbf  B$, then $\mathbf A$ and $\mathbf  B$ cannot 
be embedded in a common extension.

As another example,  let the language of $\mathcal V$
contain a constant $c$ and a unary operation $f$.
As far as the equations in $\mathcal V$
do not logically imply $f(c)=c$, 
for example, if $\mathcal V$ is axiomatized by the empty
set of equations, let $f(c)=c$
hold in $\mathbf A$ and 
$f(c) \neq c$ in $\mathbf B$.
Then $\mathbf A$ and $\mathbf  B$ cannot be embedded
in a common extension.

(b) Again in order to get the joint embedding property,
it is not enough to assume that 
the equations axiomatizing $\mathcal V$ are linear.
It is necessary to assume that the equations in 
Theorem \ref{thm} are equilinear, 
namely, that in each equation the left-hand and right-hand sides
contain exactly the same variables.

Let the language of $\mathcal V$ have two unary operations
$f$, $g$ and let $\mathcal V$ be axiomatized by the equation
$f(x)=f(y)$, which is linear but not
equilinear. For  every algebra $\mathbf A$ in
$\mathcal V$, the set $\{ a \in A \,  \mid f(a')=a,
\text{ for some } a' \in A  \,\}$ has exactly one element,
call it $a$. Let $\mathbf  B$ be another algebra in
$\mathcal V$ and let $b \in B$ be the unique element
such that $f(b')=b$, for every $b' \in B$.
If $g(a)=a$ in  $\mathbf A$   
 and $g(b) \neq b$  in $\mathbf  B$,
then obviously $\mathbf A$ and $\mathbf  B$ cannot be embedded in a common extension. The formal computation goes as follows.
Were $\iota$ and $\kappa$  embeddings
from $\mathbf A$ and $\mathbf  B$, respectively,
to some algebra $\mathbf  D$, then
\begin{align*}
  \iota(a)&=\iota(g(a))=
\iota(g(f(a))) = g(f(\iota(a))) =
g(f( \kappa (b))) = \kappa ( g(f( b))) =
\kappa ( g( b))
\\
& \neq \kappa ( b) =
  \kappa ( f(b)) =  f(\kappa ( b)) = 
f(\iota ( a))  = \iota (f(a)) = \iota (a) ,
 \end{align*}
a contradiction.

Notice that, had $\mathbf A$ and $\mathbf  B$
a common subalgebra $\mathbf  C$, then the above argument
could not be carried over, since then
$\mathbf  C$ has an element 
$c$ such that $f(c)=c$, hence  
$a=c=b$. 
Of course, the arguments in (a)
cannot be carried over, either, when 
$\mathbf A$ and $\mathbf  B$ have
a common subalgebra $\mathbf  C$, since then constants are already
interpreted in $\mathbf  C$.
 \end{remark}  

The above counterexamples are rather pathological.
Indeed, as we mentioned,
the counterexamples prevent
the joint embedding property,
but we are still able to get
the strong amalgamation property.
We first need to introduce some terminology.

\begin{definition} \labbel{exce}
Let $\mathcal V$ be a variety
and $t(z_1, \dots, z_p)$
be a term such that the variables
 $z_1, \dots, z_p$ are pairwise distinct.
More formally, the intended meaning of the
above sentence is that we write $t(z_1, \dots, z_p)$
 to mean that $t$ is a term and 
$\{ z_1, \dots, z_p\}$ is the set of variables occurring in 
$t$, enumerated without repetitions.
Subsequently, we might substitute some variables
in $t$ for other variables: the above convention applies
only to the original term (after a substitution, we formally
get a different term). 

 Under the above convention, call a variable $z_i$ \emph{exceptional (in $t$ for $\mathcal V$)}
if   
\begin{equation}\labbel{exc}
    t(z_1, \dots, z_{i-1},z_i,z_{i+1}, \dots, z_p)
= t(z_1, \dots, z_{i-1}, w, z_{i+1},\dots, z_p)
   \end{equation}
is an equation valid in $\mathcal V$, where 
$w$ is a variable distinct from  $z_i$. Otherwise, the
variable $z_i$ is \emph{ordinary}. 

Notice that it is no loss of generality to assume that 
$w$ in \eqref{exc} is distinct from
all the $z_j$'s, $j =1, \dots, p$.
Indeed, if     
$t(z_1, \dots, z_i, \dots, z_p)
= t(z_1, \dots, z_j, \dots, \allowbreak z_p)$,
with $i \neq j$, is an equation valid in $\mathcal V$,
then, applying twice the equation, we get
\begin{equation} \labbel{ban}      
t(z_1, \dots, z_i, \dots, z_p)
= t(z_1, \dots, z_j, \dots, z_p)=
t(z_1, \dots, w, \dots, z_p)
 \end{equation}
since $z_i$ occurs only once
on the left-hand side and, symmetrically,
$w$ occurs only once
on the right-hand side, since we have assumed that 
$w$  is distinct from
all the $z_j$'s. 
Here we have used the assumption
that all the  variables in the original term  $t$ 
on the left-hand side 
are distinct.
The above notions are dependent on the variety $\mathcal V$,
but the variety $\mathcal V$ shall be always clear from the context,
henceforth we shall drop reference to $\mathcal V$.

Of course, it might happen that some term has no
exceptional variable and, less frequently,
that it has no
ordinary variable.
\end{definition}

 \begin{definition} \labbel{lingen}    
In our general definition of a linear equation
we require that there is at most one occurrence
of operations on each side.
If the language under consideration has possibly constants
and since we do not consider constants as operations,
the general forms of a \emph{linear} equation 
are as follows.
\begin{align} \labbel{linec}  
f(w_{1},w_{2}, \dots, w_{n} ) &= v_1,  \text{ or } 
\\
\labbel{dc}  
f(w_{1},w_{2}, \dots, w_{n} ) &= g(v_{1},v_{2}, \dots, v_{m}    ),
 \end{align}
where $f$ and $g$ are
 operations
and each one among $w_i$ or $v_j$  is either a variable or a constant. 
\end{definition}

Strictly speaking, also the equation  $w_1=v_1$
obeys our definition of a linear equation.
However, if either $w_1$ or  $v_1$ is a variable,
then we are in a trivial variety, in which all the statements we consider 
are either trivially true or trivially false.  
If both $w_1$ and  $v_1$ 
are constants, we can equivalently consider a variety 
without the constant $v_1$, substituting every occurrence of $v_1$
for $w_1$  in every equation. In this respect, compare also a remark
at the beginning of the proof of the following theorem.

\begin{theorem} \labbel{thmc}
If $\mathcal V$ is a variety 
axiomatized by a set of linear equations,
possibly with constants, 
then $\mathcal V$ has the strong amalgamation property.
 \end{theorem}

 \begin{proof}
The proof is similar to the proof of Theorem \ref{thm},
but in the present case there are more situations
in which the value of some 
expression of the form
$f(d_1, d_2, \dots )$ is ``forced''
by some condition. 

In detail, 
let $\mathbf A$, $\mathbf  B$ and $\mathbf  C$
be algebras to be amalgamated as in \eqref{apeq}.
We shall highlight the various steps
of the proof.

\emph{A preliminary remark about constants.} 
First of all, notice that 
if $c$ is a constant, then $c$
is interpreted by the same element in 
$\mathbf A$, $\mathbf  B$ and $\mathbf  C$,
since  $\mathbf  C$ is a subalgebra of 
both $\mathbf A$ and $\mathbf  B$.
Hence it is not necessary to explicitly state where 
$c$ is interpreted and, as a usual abuse of notation,
we might not graphically distinguish 
the constant symbol from its interpretation,
though we shall actually do it in a paragraph below
for the sake of clarity.  

Moreover, if two constants
$c_1$ and $c_2$  are interpreted by the same element,
say, in $\mathbf A$, then 
$c_1$ and $c_2$  are interpreted by the same element
in $\mathbf  C$, hence in $\mathbf  B$, too.
Thus $c_1$ and $c_2$  should be interpreted by the same element
in any amalgamating algebra $\mathbf  D$, as well. 
Hence, in the above situation, it is no loss of generality to
assume that $\mathcal V$ satisfies the equation
$c_1 = c_2$. Again without loss of generality,
we can work in the variety obtained from
$\mathcal V$ by deleting the constant $c_2$
from the language and replacing $c_2$ by $c_1$
in every equation axiomatizing $\mathcal V$.   Thus if we prove the theorem
in the special case of
 algebras in which no two distinct constants are interpreted by the same
element, we get  a proof for the full theorem.
Henceforth, we shall assume that 
no pair of constants are interpreted by the same
element.

\emph{A term associated to an expression $f(d_1, d_2, \dots , d_n)$.} 
We are now going to address
the following problem. 
Given
an expression  $f(d_1, d_2, \dots , d_n)$,
with $f$ an operation in the language of $\mathcal V$ 
and  $d_1, d_2, \dots, d_n \in A \cup B$,
when is the value of 
$f(d_1, d_2, \dots, d_n )$ forced by some condition?
Of course, the first three clauses 
in \eqref{alf} in the proof of Theorem \ref{thm}
force the value of  $f(d_1, d_2, \dots, d_n )$,
but when considering linear not
equilinear equations there are more possibilities.
For example, if
$f(x,x,y)= f(x,x,z)$ holds in 
$\mathcal V$, 
then   the first three clauses 
in \eqref{alf} do not necessarily  force the value of
$f(a,a,b)$, when $a \in A \setminus B$
and $b \in B \setminus A$;
   however, the equation 
$f(x,x,y)= f(x,x,z)$
implies 
$f(a,a,b)=f(a,a,a)$
and the latter expression can be evaluated
in $\mathbf A$.
We  now  present
the general form of the argument. 

Consider again the  expression 
$f(d_1, d_2, \dots, d_n )$, with
$f$  an operation in the language of $\mathcal V$
and $d_1, d_2, \dots, d_n \in A \cup B$. 
In the expression
$f(d_1, d_2, \dots, d_n )$
replace every element interpreting some constant
by the corresponding constant symbol, 
and replace all the other elements by variables
in a bijective way, namely, to distinct elements
make correspond distinct variables
 and to equal elements
make correspond equal variables.
Consider the expression obtained in this way
as a term in which no variable is repeated.
We thus get 

(*)  a term  $t(x_{1}, x_{2}, \dots)$ such that 
  $x_{1}, x_{2}, \dots$
are pairwise distinct variables
and, in any amalgamating algebra,
 $f(d_1, d_2, \dots )$
is expressed as 
$t(d_{i_1}, d_{i_2}, \dots)$,
where $d_{i_1}, d_{i_2}, \dots$ 
are pairwise distinct elements
from the set 
$\{ d_1, d_2, \dots\}$ and 
no $d_{i_r}$ is the interpretation of some constant. 

The intended meaning in the construction of $t$
is that  we should have $f(d_1, d_2, \dots )=t(d_{i_1}, d_{i_2}, \dots)$ but,
of course, neither $f(d_1, d_2, \dots ) $ nor $ t(d_{i_1}, d_{i_2}, \allowbreak  \dots)$
are generally interpretable, until we have constructed some amalgamating
algebra. On the other hand, we do have 
$f(d_1, d_2, \dots )=t(d_{i_1}, d_{i_2}, \dots)$ in $\mathbf A$,
in case all the elements $d_1, d_2, \dots $ belong to $A$,
and similarly for $\mathbf  B$.

More formally, the above term
$t$ associated to the expression
$f(d_1, d_2, \dots, \allowbreak d_n )$ is introduced as follows.
 For every  expression of the form  
$f(d_1, d_2, \dots, d_n )$, let
$I=\{ 1,2,\dots, n\}$ and 
\begin{equation*}\labbel{iiii}    
I'= \{ \, i  \in I \mid d_i \text{ is 
not the interpretation of some constant} \, \}.
  \end{equation*} 
Define an equivalence relation $ \sim'$ on $I'$ by   
\begin{equation*}\labbel{puhhh}
i \sim' j \quad \text{ if } \quad d_i=d_j. 
  \end{equation*}    
To distinct equivalence classes of $ \sim'$
associate  distinct variables.
Then $t$ is defined by
\begin{equation}\labbel{pupu}
t(x_{1}, x_{2}, \dots)= f(w_1, w_2, \dots, w_n),
  \end{equation}     
where, for $i=1, \dots, n$, 
if $d_i$ is the interpretation of some constant, then  $w_i$
is the  symbol for that constant; otherwise 
$w_i$ is the variable associated to the $ \sim'$-class
 of $i$.

The term $t$ is essentially unique, modulo
the choice or rearrangements 
 of variables. Notice that, by the  assumption 
in the preliminary remark above, 
every element can be the interpretation of at most
one constant, hence $t$ is well-defined. 
For example, given the expression
$f(d,d,e,e,g,g,g)$,
where $g=c_ { {_\mathbf  C}}$ is 
the interpretation of some constant $c$, 
we  take
$t(x,y) $ to be the term $ f(x,x,y,y,c,c,c)$,
thus 
$f(d,d,e,e,g,g,g)=
f(d,d,e,e,c_ { {_\mathbf  C}},c_ { {_\mathbf  C}},c_ { {_\mathbf  C}})=t(d,e)$
in any amalgamating algebra.
Had we  taken
$t'(x,y,z) $ as $ f(x,x,y,y,z,z,z)$,
we still have
$f(d,d,e,e,g,g, \allowbreak  g)
=t'(d,e,g)$, but in this case
we do not meet the
 requirement that $g $
is not the interpretation of some constant.
Notice also that we want  the elements
$d_{i_1}, d_{i_2}, \dots$
appearing in the (possible) evaluation of $t$ to be pairwise distinct.
Hence, in the above situation,
defining  
$t(x,x',y) = f(x,x',y,y,c,c,c)$
again would not meet our requirements, since then,
in order to get $f(d,d,e,e,c,c,c)$,
we need to express it as
$t(d,d,e)$.  

\emph{When is the value of $f(d_1, d_2, \dots )$
forced by some condition?} 
Now we can answer the above question. 
For every operation $f$,
and every sequence  
$ d_1, d_2, \dots \in A \cup B$
of appropriate length,
we say that 
the value of $f(d_1, d_2, \dots )$
is \emph{forced} if one of the
following two eventualities occurs.
  \begin{enumerate}[(E1)]    
\item 
For the term  $t$  constructed above in (*)
depending on $f$ and $ d_1, d_2, \dots$,
the equation  $t(x_{1}, x_{2}, \dots)=v$
is  valid in $\mathcal V$,
for some
 $v$ either a constant $c$ or some variable $x_j$ 
among $x_{1}, x_{2}, \dots$. As in the proof of \ref{thm},
if such a $v$ exists, it is unique, unless we are in a trivial variety, 
in which case the theorem is obviously true.
Notice that we have excluded in the preliminary remark the possibility 
that $c_1=c_2$ for two distinct constants, hence we cannot have
at the same time $t(x_{1}, x_{2}, \dots)=c_1$ and
$t(x_{1}, x_{2}, \dots)=c_2$.
Excluding the  trivial variety, 
under the present eventuality
the \emph{value forced for} $f(d_1, d_2, \dots )$
is  either the interpretation
of the constant $c$, or  the element from $\{  d_1, d_2, \dots \}$ 
the variable    $x_j$ has been substituted for in the construction
of $t$.
\item 
Let $t$  be the term from  (*)
and suppose, without loss of generality, 
that the variables 
are rearranged  in such a way 
that all the ordinary variables appear before
the exceptional variables,
as in  $t(x_{1}, x_{2}, \dots, x_{ \ell -1}; \allowbreak 
 x_{ \ell}, \dots, x_r)$, where the semicolon 
separates the two kinds of variables.
Let $d_{i_1}, d_{i_2}, \dots,
d_{i_{ \ell -1}}$ be the elements corresponding to the ordinary
variables $x_{1}, x_{2}, \dots, x_{ \ell -1}$.
In this case the value of 
 $f(d_1, d_2, \dots )$ is \emph{forced} if 
either $\{ d_{i_1}, d_{i_2}, \dots,
d_{i_{ \ell -1}} \} \subseteq A$
or  $\{ d_{i_1}, d_{i_2}, \dots,
d_{i_{ \ell -1}} \} \subseteq B$.
If one of these two eventualities occurs,
the \emph{value forced for} $f(d_1, d_2, \dots )$
is
$t(d_{i_1}, d_{i_2}, \dots,
d_{i_{ \ell -1}};
 d, \dots, d)$, evaluated, respectively, in $\mathbf A$ 
or in $\mathbf  B$, where $d$ is an arbitrary element
of $C$.  
  \end{enumerate}

A few remarks are in order about the latter clause.
By applying several times
the definition of an exceptional variable, we get that the equations
\begin{equation*}     
t(x_{1},  {\dots}, x_{ \ell -1};
 x_{ \ell}, {\dots}, x_p)= 
 t(x_{1},  {\dots}, x_{ \ell -1};
 z, {\dots}, z) = t(x_{1}, {\dots}, x_{ \ell -1};
 z', {\dots}, z')
  \end{equation*}
 hold in $\mathcal V$, hence
the definition of the value  forced
does not depend on the choice of $d$.  
Notice that it is important to have $C$
nonempty, in order for this value to be well-defined.
We might be in the extreme case when all the variables
in $t$ are exceptional; then, without a common subalgebra $\mathbf  C$
 of both $\mathbf A$ and $\mathbf  B$, the value of 
$t(a,a, \dots )$ computed in $\mathbf A$ might be different
from the value of   $t(b,b, \dots )$ computed in $\mathbf  B$.
Compare Remark \ref{con}(b). 
Of  course, it may happen that 
$\{ d_{i_1}, d_{i_2}, \dots,
d_{i_{ \ell -1}} \} \subseteq C$ in (E2).
If this is the case, then
$t(d_{i_1}, d_{i_2}, \dots,
d_{i_{ \ell -1}};
 d, \dots, d)$ is evaluated in the same way in $\mathbf A$, 
$\mathbf  B$ and $\mathbf  C$, since $\mathbf  C$
is a subalgebra of both $\mathbf A$ and $\mathbf  B$.
Hence the definition is not ambiguous, in this case, too.
Finally, notice that the first two clauses 
in \eqref{alf} from the proof of Theorem \ref{thm}
 are actually a special case of the above eventuality (E2). 

We now check that (E1) and (E2) force the same value for
$f(d_1, d_2, \dots )$ in any overlapping case.
Indeed, if  (E1)  applies,
then   
$v=t(x_{1},  \dots, x_{ \ell -1};
 x_{ \ell}, \dots, \allowbreak  x_p) = 
 t(x_{1},  \dots, x_{ \ell -1};
 z, {\dots}, z)$
 are equations valid in $\mathcal V$.
If $v$ is a constant $c$ and also (E2) 
applies, then
the above equations imply  
$c=t(d_{i_1}, d_{i_2}, \dots,
d_{i_{ \ell -1}};
 d, \dots, \allowbreak d)$, as evaluated
either  in $\mathbf A$ or in $\mathbf  B$.
The argument is similar if
$v$ is variable among  
the ordinary variables $x_{1},  {\dots}, x_{ \ell -1}$ of $t$.
Otherwise, $v$  is an exceptional variable, but then 
$v=t(x_{1},  {\dots}, x_{ \ell -1};
 x_{ \ell}, {\dots}, x_p)$  together with \eqref{ban}
imply that we are in a trivial variety. 

\emph{Defining the amalgamating algebra.} 
Since all the above ``forced'' conditions
follow from equations valid in $\mathcal V$,
it is evident that in any amalgamating algebra $\mathbf  D$ 
the expression $f(d_1, d_2, \dots )$ 
should  actually assume the forced  value indicated
in (E1) or (E2), if applicable. 
Now, as in the proof of \ref{thm},
we pick an arbitrary element 
$ \bar{d} \in A \cup B$ 
and let 
$f(d_1, d_2, \dots )= \bar{d}$
if no value is forced by either (E1) or (E2).
In detail, we let $D=A \cup B$ and,
for $d_1, d_2, \dots \in D$
and $f$ an operation in the language of $\mathcal V$,   we define
 $f(d_1, d_2, \dots )$ by
\begin{equation}\labbel{alfc}    
 \begin{aligned} 
  f(d_1, d_2, \dots ) &= \text{the value forced by (E1) or (E2),}
\\[-0.2 pt]
&\phantom{ = }\text{ \ if either (E1) or (E2) applies, and} 
\\
 f(d_1, d_2, \dots ) &= \bar{d},  
\text{\qquad\qquad\qquad otherwise,} 
   \end{aligned}
\end{equation}   
thus obtaining an algebra $\mathbf  D$ appropriate for the language
of $\mathcal V$.

As we mentioned, (E1) and (E2)
force the same value in any overlapping case,
hence the definition is correct.
Moreover, when 
$d_1, d_2, \dots \in A$, eventuality (E2)
applies and gives the same value for  
 $f(d_1, d_2, \dots )$ as evaluated in $\mathbf A$,
by a comment shortly after the definition (*) of $t$.  
Hence $\mathbf  D$ extends $\mathbf A$.
Similarly, $\mathbf  D$ extends $\mathbf B$.

\emph{Checking that the equations are satisfied. Case \eqref{linec}.} 
It remains to show that 
the equations of the form 
\eqref{linec} and \eqref{dc} 
 axiomatizing $\mathcal V$ hold in $\mathbf  D$.
Eventuality (E1) takes care of equations 
of the form 
\eqref{linec}.
We first outline the argument.
Suppose that 
some equation of the form \eqref{linec}
is evaluated in $\mathbf  D$ modulo some assignment.
Of course, when we evaluate
$f(w_{1},w_{2}, \dots, w_{n} )$
modulo an assignment  
$y_1 \mapsto d_1$,   
$y_2 \mapsto d_2$,\dots,
where $y_1, y_2, \dots$
are the variables (enumerated without repetitions) in the set $\{w_{1},w_{2}, \dots, w_{n} \}$,  
it might happen that distinct variables
are assigned equal elements;
moreover, some variable might be assigned
the value of some constant, hence the term 
$t$ from (*) might contain a smaller number of variables, say,
 it has the form $t(x_{1}, x_{2}, \dots)$,
where $ |\{ x_{1}, x_{2}, \dots \}| \leq |\{ y_{1}, y_{2}, \dots \}|$.
However, if we substitute all the occurrences of some variable
by the same constant symbol in \eqref{linec}, we still get an equation
valid in $\mathcal V$. Similarly, we can substitute many variables
for a single different variable in  
\eqref{linec}, getting an equation valid in $\mathcal V$.
Thus, under substitutions as above,
$f(w'_{1},w'_{2}, \dots, w'_{n} ) = v'_1$
is an equation valid in $\mathcal V$.

More formally,
suppose that 
$f(d_{1},d_{2}, \dots, d_{n} )$
is the evaluation of 
$f(w_{1},w_{2}, \allowbreak \dots, w_{n} )$
under the assignment $\rho$ given by 
$y_1 \mapsto d_1$,   
$y_2 \mapsto d_2$,\dots,
where $y_1, y_2, \dots$
are the variables in the set $\{w_{1},w_{2}, \dots, w_{n} \}$.
If \eqref{linec} $f(w_{1},w_{2}, \dots, \allowbreak  w_{n} )= v_1$
is an equation valid in $\mathcal V$, 
we have to show that in $\mathbf  D$ 
$f(d_{1},d_{2}, \dots, d_{n} )  \allowbreak = \rho (v_1)$
if $v_1$ is a variable, and that  
$f(d_{1},d_{2}, \dots, d_{n} )=c$
if $v_1$ is the constant $c$. 
Let $\tau$ be the function from
$\{ d_1, d_2, \dots, d_n\}$
which has been implicitly used in the definition 
of the term $t$, namely,  
if $d_i$ is the interpretation of some constant $c$, then 
$\tau(d_i)$ 
is the  symbol for that constant; otherwise 
$\tau(d_i)$ is the variable associated to the $ \sim'$-class
 of $i$.  Then \eqref{pupu} reads  
\begin{equation}\labbel{pupux}
t(x_{1}, x_{2}, \dots)= f( \tau (d_1), \tau ( d_2), \dots, \tau (d_n)).
  \end{equation}     
Extend $\rho$ by sending a constant symbol to its interpretation.
Since $f(w_{1},w_{2}, \allowbreak   \dots,   w_{n} )= v_1$
is  valid in $\mathcal V$, then also 
\begin{equation}\labbel{pupuxx}
f( \tau ( \rho (w_1)), \tau ( \rho ( w_2)), \dots, \tau ( \rho (w_n)))
=  \tau ( \rho (v_1))
  \end{equation}     
is valid in $\mathcal V$.
By \eqref{pupux} and \eqref{pupuxx}, the equation 
\begin{equation}\labbel{pupups}
t(x_{1}, x_{2}, \dots)= \tau ( \rho (v_1))
  \end{equation}     
is valid, as well, hence 
we can apply eventuality (E1)
getting 
$f(d_{1},d_{2}, \dots, d_{n} )= \rho (v_1)$,
what we had to show. 

\emph{Checking that the equations are satisfied. Case \eqref{dc}.} 
It remains to consider equations
of the form \eqref{dc}.   
We need to prove that 
$f(d_{1},d_{2}, \dots  ) = g(e_{1},e_{2}, \dots  )$,
when $d_{1},d_{2}, \dots, e_{1},e_{2},\dots \in A \cup B$
and equal elements are substituted for equal variables
in some equation from \eqref{dc} valid in $\mathcal V$.
We first sketch the argument.
Associate to $f(d_{1},d_{2}, \dots  )$ 
and  $g(e_{1},e_{2}, \dots  )$
corresponding terms $t$ and  $s$ as in (*),
 naming the variables in a
 consistent way, say, if $d_j=e_k$,
then the variable corresponding to  
$d_j $ in the construction of $t$
should be equal to the variable 
corresponding to  
$ e_k$ in the construction of $s$.
As in the previous paragraph, the validity 
of \eqref{dc} in $\mathcal V$ implies the validity of 
an equation of the form
$t=s$.
Thus the value of 
$f(d_{1},d_{2}, \dots  ) $ is forced by
(E1) if and only if 
the value of  $  g(e_{1},e_{2}, \dots  )$
is forced by (E1).
If we check that
$t$ and  $s$ have the same ordinary variables, then   
the corresponding statement
holds with regards to (E2), 
hence 
the value of 
$f(d_{1},d_{2}, \dots  ) $ is forced if and only if 
the value of  $  g(e_{1},e_{2}, \dots  )$
is forced.
If neither value is forced, then both expressions
are assigned the value $ \bar{d}$. In all cases, they are equal.

In order to write explicitly the details, it is convenient to write 
\eqref{dc} as   
\begin{equation}\labbel{dcx}
f(w_{1},w_{2}, \dots, w_{n} ) = g(w_{n+1},w_{n+2}, \dots, w_{p}),
   \end{equation}    
where each $w_i$   is either a variable or a constant and, similarly, 
express the identity to be proved as 
\begin{equation}\labbel{dubba}    
f(d_{1},d_{2}, \dots,  d_{n}) = g(d_{n+1},d_{n+2}, \dots,  d_{p}).
  \end{equation}
As in the previous case, let  \eqref{dubba}  be obtained 
from \eqref{dcx} through some assignment  $\sigma$.
As in the definition (*) of the associated term, 
let
$I=\{ 1,2,\dots, p\}$, 
$ I'= \{ \, i  \in I \mid d_i \text{ is 
not the interpretation of some constant} \, \}$
and let  $ \sim'$ on $I'$ be defined by   
$i \sim' j $ if $d_i=d_j$. 
To distinct equivalence classes of $ \sim'$
associate  distinct variables, say, 
$\tau(d_i)$ is the variable associated 
to the class of $d_i$. If $i \in I \setminus I' $,
let    $\tau(d_i)$ be the constant of which $d_i$  
is the interpretation.
Thus, as in \eqref{pupux},  the term  $t$ associated to the expression
$f(d_{1},d_{2}, \dots,  d_{n}) $ 
is 
$t(x_{1}, x_{2}, \dots)= f( \tau (d_1), \tau ( d_2), \dots, \tau (d_n))$
and, similarly, 
the term  $s$ associated to the expression
$  g(d_{n+1},d_{n+2}, \dots,  d_{p})$ is
$s(y_{1}, y_{2}, \dots)= g( \tau (d_{n+1}), \tau ( d_{n+2}), \dots, \tau (d_{p}))$.
Of course, it might happen that the sets 
$\{ x_{1}, x_{2}, \dots \}$ and  
$\{ y_{1}, y_{2}, \dots \}$ are distinct; however, we 
shall see that their symmetrical difference consists 
entirely of exceptional variables. As we mentioned,
the choice of variables for $s$ and  $t$ here is not arbitrary,
we have named variables in such a way that  
if $d_i=d_j$, then the variables associated 
to $d_i $ and $ d_j$ are the same, even when
$d_i $ is in the range of $f$ and  $ d_j$
is in the range of $g$.  
Henceforth, 
\begin{equation}\labbel{idid}
t(x_{1}, x_{2}, \dots ) = s(y_{1}, y_{2}, \dots)
  \end{equation}    
is an equation valid in $\mathcal V$, since 
\eqref{dcx} is valid in $\mathcal V$.
As in the previous case, \eqref{idid}
and \eqref{dcx} are not necessarily equivalent,
since some variables from \eqref{dcx}
might have been identified in \eqref{idid}; and moreover, some variables
from \eqref{dcx} might have been replaced by constants.
However, \eqref{dcx} does imply \eqref{idid}, and this is
what we need.  

From \eqref{idid} we immediately get
that  eventuality (E1) holds for 
$f(d_{1},d_{2}, \dots, \allowbreak   d_{n}) $ 
if and only if it holds for
$  g(d_{n+1},d_{n+2}, \dots,  d_{p})$  
 and, if this is the case, the forced vales are equal.

To deal with (E2), we first need to observe that
if some variable $x$ occurs
in $t$ but not in $s$, then 
$x$ is exceptional in $t$ for $\mathcal V$,
by applying     
\eqref{idid} twice as  follows
$t(x_{1}, x_{2}, \dots, x, \dots) = s(y_{1}, y_{2}, \dots)
= t(x_{1}, x_{2}, \dots, z, \dots)$.
Symmetrically, a variable occurring in $s$
and not in $t$ is exceptional in $s$.  
 It follows obviously from \eqref{idid}
that if some variable $x$ occurs both in $t$
and in $s$, then $x$ is exceptional 
in $t$ if and only if $x$ is exceptional in  $s$.   
In conclusion, the sets of variables of $t$ and $s$
might be distinct, but the sets of \emph{ordinary} variables in  
$t$ and  $s$ are equal. 
In particular, the set 
$\{ d_{i_1}, d_{i_2}, \dots\}$ of elements corresponding to the ordinary
variables is the same both for 
$f(d_{1},d_{2}, \dots,  d_{n}) $ 
and  for
$  g(d_{n+1},d_{n+2}, \dots,  d_{p})$,
thus the value of the former is forced according to (E2)
if and only if the value of the latter is forced.
If this is the case, the two forced values are the same, again
by \eqref{idid}, which holds in $\mathcal V$, hence both in $\mathbf A$
and in $\mathbf  B$.

In conclusion, we have showed that
the value of $f(d_{1},d_{2}, \dots,  d_{n}) $ 
is forced if and only if the value of
$  g(d_{n+1},d_{n+2}, \dots,  d_{p})$
 is forced and, if this is the case, the two values are equal.
Otherwise,
$f(d_{1},d_{2}, \dots,  d_{n}) =
\bar{d} =   g(d_{n+1},d_{n+2}, \dots, \allowbreak   d_{p})$.
We have proved equality in each case, hence
$\mathbf  D$ is indeed an algebra in $\mathcal V$, thus 
a desired amalgamating algebra.
\end{proof}

\begin{remark} \labbel{fullc}
(a) Remark \ref{varr} applies with no essential modification
to the context in Theorem \ref{thmc}.  

As far as part (c) in Remark \ref{varr}
is concerned, 
let $D \supseteq  A \cup B$
and $\eta $ be a function from $  D^{fin} $ to $  D$.
If $d_1, d_2, \dots, d_n \in D$ 
and the value of $f(d_1, d_2, \dots, d_n)$
is neither  forced  by (E1) nor by (E2),
set 
\begin{equation}\labbel{for}     
f(d_1, d_2, \dots, d_n)= \eta (\{d_{i_1}, d_{i_2}, \dots \})
  \end{equation}
in $\mathbf  D$, where
$\{d_{i_1}, d_{i_2},  \dots \}$ corresponds to  the set of 
ordinary  variables in the term associated to 
the expression $f(d_1, d_2, \dots, d_n)$. 
We have showed in the proof of 
Theorem \ref{thmc} that if some equation of the form
\eqref{dcx} is valid in $\mathcal V$, then, for any instance
of \eqref{dubba}, the sets of ordinary variables 
associated to $f(d_1, d_2, \dots, d_n)$
and to $g(d_{n+1},d_{n+2}, \dots,  d_{p})$
coincide, hence,
if  these values are not forced, we get
 $f(d_{1},d_{2}, \dots  )=  \eta (\{d_{i_1}, d_{i_2}, \dots \})
=g(d_{n+1},d_{n+2}, \dots,  d_{p})$ from
\eqref{for}.

(b) The full version of Corollary \ref{v}
follows from Theorem \ref{thmc}. 
Only (1) needs a comment, since we cannot use JEP to prove it.
However Theorem \ref{thmc} 
can be applied, together with the methods indicated in Remark \ref{varr}(e).
 \end{remark}   

Remark \ref{con} essentially provides 
all kinds of counterexamples to  JEP. We obtain JEP
 when the equations axiomatizing some variety 
are linear, provided there is exactly one constant
and this constant is always interpreted as a one-element
subalgebra.

\begin{corollary} \labbel{corthmc}
Suppose that $\mathcal V$ is a variety 
axiomatized by a set of linear equations, that
the language of $\mathcal V$ contains exactly one constant $c$
and that $f(c,c,c, \dots )=c$ is valid in $\mathcal V$,
for every operation $f$ in the language.  

Then $\mathcal V$ has the joint embedding property.
The class of finite algebras in $\mathcal V$ has the joint embedding
property. If the language of $\mathcal V$  is finite
(or just if there are countably many finite algebras up to isomorphism),
then  the class of finite algebras in $\mathcal V$ has  a Fra\"\i ss\'e limit.
 \end{corollary}

\begin{proof}
Under the assumptions, 
the one-element algebra is (isomorphic to) a subalgebra of every algebra.
Hence AP implies JEP, by taking $\mathbf  C$ a one-element algebra. 
\end{proof}

The above results imply that many varieties 
corresponding to Maltsev conditions ``localized''
at some constant $0$ have SAP, JEP and Fra\"\i ss\'e limits.
In particular, Theorem \ref{thmc} and 
Corollary \ref{corthmc} apply to varieties 
with operations witnessing permutability at $0$,
arithmeticity at $0$
\cite[Theorems 1(2) and 4(2)]{Du}, distributivity at $0$
\cite[Theorem 1]{Ch1}, 
$3$-permutability at $0$ \cite[Theorem 1]{Ch2},
 $n$-permutability at $0$ together with distributivity at $0$
\cite[Theorem 2(3)]{Ch2}, as well as varieties 
with a near $0$-unanimity operation 
\cite[Definition 4]{Du}. 
In this respect, compare also the example after Definition 
 11.1 in \cite{Gu}.

\begin{remark} \labbel{quas}
We cannot generalize Theorem \ref{thmc}
to \emph{quasiequations} (= Horn sentences in a language without relations),
 namely, if some quasivariety $\mathcal Q$ 
is axiomatized by means of quasiidentities constructed using linear identities,
it is not necessarily the case that $\mathcal Q$ has AP.

For example, let $\mathcal Q$ be the quasivariety in the language with two
unary operations $f$ and  $g$ such that  $\mathcal Q$ is 
axiomatized by the quasiequation    
 \begin{equation}\labbel{q}
f(x)=f(y) \quad \Rightarrow \quad  x=y.
  \end{equation}     
Then $\mathcal Q$ has not AP.

Indeed, let $\mathbf  C$ be $\mathbb N^+$
with both $f$ and  $g$ interpreted as the successor function.
Extend $\mathbf  C$ to $\mathbf A$ and $\mathbf  B$ by adding
two copies of $0$. In $\mathbf A$ we let   
$f$ and  $g$ be again the successor function.
On the other hand, in $\mathbf  B$,
we let $f$ be the successor function, while we 
set $g(0')= 0'$.
If by contradiction there is an amalgamating structure, we  have
$f(0)=f(0')$, hence $0$ and  $0'$ should be identified by \eqref{q}, but this is 
impossible, since   $g(0) \neq g(0')$.

Of course, it might happen that some analogue 
of Theorem \ref{thmc} holds for 
quasiequations of some special form, but we have not
investigated in depth this possibility. 
\end{remark}

\section{Adding unary operations} \labbel{add} 

It is easy to see that if some variety has the amalgamation property,
then the amalgamation property is preserved by adding
a new unary operation $h$ with equations asserting that 
$h$ is an endomorphism. See the following observation.

In the next observation $\mathcal V$ is an arbitrary variety
with the amalgamation property. We are not assuming that $\mathcal V$ 
is axiomatized by linear equations.
The observation is probably folklore, but we know no reference
for it.

\begin{observation} \labbel{obs}
Suppose that $\mathcal V$ is a variety with the (strong) amalgamation
property. Let $\mathcal V_h$
be the variety in a language with a new unary operation
$h$ added; the equations axiomatizing   
 $\mathcal V_h$ are the equations axiomatizing $\mathcal V$ 
plus equations asserting that $h$ is an endomorphism of $\mathcal V$, namely,
\begin{equation}\labbel{pup}
h(f(x_1,x_2, \dots ))=f(h(x_1),h(x_2), \dots ) 
   \end{equation}    
for every operation $f$ in the language of 
$\mathcal V$ and where $x_1,x_2, \dots $
are distinct variables.  

Then $\mathcal V_f$ has 
the (strong) amalgamation
property. 
 \end{observation}

\begin{proof}
The argument is purely categorical.
If $\mathbf A, \mathbf  B, \mathbf  C \in \mathcal V$ are 
as in the assumptions of (S)AP and have an amalgam in $\mathcal V$, then
their push out $\mathbf  D$  (in the category of algebras in $\mathcal V$ with
 $\mathcal V$-homomorphisms) is still an amalgam \cite{KMPT}.

Now suppose that $\mathbf A_h, \mathbf  B_h, \mathbf  C_h $  are
in $  \mathcal V_h$, let 
$\mathbf A, \mathbf  B, \mathbf  C $ be their reducts
to the language of $  \mathcal V$
and let $\mathbf  D$ be as above.
The unary operation $h$ satisfies \eqref{pup} on 
$\mathbf A_h $, $\mathbf  B_h$ and $\mathbf  C_h $, hence
we can think of $h_{\mathbf A}$, $h_{\mathbf B}$ and
$h_{\mathbf C}$
 as  $\mathcal V$-endomorphisms of
$\mathbf A, \mathbf  B, \mathbf  C $, respectively.
Composing these endomorphisms with the appropriate embeddings,
we get a commuting diagram of homomorphisms from $\mathbf A$,
$\mathbf  B$ and $\mathbf  C$  to $\mathbf  D$.
By the push out property of $\mathbf  D$,
 such homomorphisms  can be extended to an endomorphism 
$h_{\mathbf D}$
of $\mathbf  D$.
Now turn back and think of 
$h_{\mathbf D}$ as a
unary operation on  $D$.
Expanding $\mathbf  D$ by adding 
$h_{\mathbf D}$ we get an algebra in 
$\mathcal V_f$ which amalgamates  
$\mathbf A_h, \mathbf  B_h, \mathbf  C_h $.
\end{proof}
 
\begin{remark} \labbel{rmo}   
(a) Of course, Observation \ref{obs}
holds in the case we add a family of
unary operations which are 
 endomorphisms. 

(b) Observation \ref{obs} applies also to an arbitrary class $\mathcal K$ 
of structures (not necessarily a variety),
provided $\mathcal K$ has push outs.
In the general case, condition \eqref{pup}
 should be replaced by a condition asserting that $h$
is an endomorphism. 

(c) We can also add unary operations which
are automorphisms.
The property that some $h$ is bijective
cannot be expressed equationally without expanding the language;
however, we might assume that 
there is one more unary operation
$k$ and that 
$k(h(x))= h(k(x))=x$
are equations  valid in the variety under consideration.
\end{remark}   

The assumption that \eqref{pup}
holds  
for \emph{every} operation $f$ in the language of 
$\mathcal V$ is necessary in Observation
\ref{obs}. Let $+$ and juxtaposition denote lattice
operations. 

\begin{example} \labbel{ex}
The variety of distributive lattices with a unary operation
$h$ satisfying 
\begin{equation}\labbel{latt}
h(x+y)= h(x)+h(y)
  \end{equation}
has not AP.     
 \end{example}   

\begin{proof}
Let $\mathbf  C$ be the three-elements
chain with base set $\{ 0, 1, 2\}$,
and let $h$ be the identity on $C$.

Expand $\mathbf  C$ in two ways.
Let $A= \{ 0, 1, 2, a\}$,
with $a$ a complement of $1$ in $\mathbf A$
and, again, let $h$ be the identity on $A$.   
Let $B= \{ 0, 1, 2, b\}$,
with $b$ a complement of $1$ in $\mathbf B$
and let $h$ be the identity on $C$
and $h(b)=2$.

Then $\mathbf  C \subseteq \mathbf A, \mathbf  B$
and all algebras satisfy \eqref{latt}.
However, in any amalgamating structure, 
$a$ and $b$ should be identified, since complements are unique
in distributive lattices.   
This contradicts $h(a)=a$ and  $h(b)=2$,
hence $\mathbf A$, $\mathbf  B$ and $\mathbf  C$ cannot
be amalgamated.
\end{proof}

The variety of distributive lattices has the amalgamation
property, but not the strong amalgamation property \cite{FG}.
In fact, Example \ref{ex} is just a small elaboration on this result. 

\begin{problem} \labbel{prob}
Is there an example similar to Example \ref{ex}
in which $\mathcal V$ has the strong amalgamation property?

Is there an example similar to Example \ref{ex}
in which $h$ is bijective? (We can express
the assumption that $h$ is bijective by introducing a further operation
$k$, as in Remark \ref{rmo}(c).)
\end{problem}  

We have seen in Example \ref{ex}
that Observation \ref{obs}
does not generalize when we assume that
$h$ preserves only some operation in $\mathcal V$. 
However, the observation does generalize when 
$\mathcal V$ is axiomatized by linear equations
and $h$ is assumed to be bijective.

\begin{proposition} \labbel{un}
Suppose that $\mathcal V$
is a variety axiomatized by a set of linear equations and let 
$\mathcal V_{h,k}$ be a variety with two new unary operations
$h$ and  $k$ added.   
The variety $\mathcal V_{h,k}$ is supposed to satisfy all the
equations of $\mathcal V$, the equations
\begin{equation} \labbel{kh}  
k(h(x))= h(k(x))=x, \qquad h(c)=c,
  \end{equation}    
for every constant $c$ in the language of $\mathcal V$, 
as well as possibly some equations of the form \eqref{pup},
$f$ varying on a set of operations of $\mathcal V$.  

Then $\mathcal V_{h,k}$ has the strong amalgamation property.

If either $\mathcal V$ is 
axiomatized by a set of equilinear equations
in a language without constants, or  the language of 
$\mathcal V$ has only one constant which represents a subalgebra,
then $\mathcal V_{h,k}$ has the joint embedding property.
 \end{proposition}

  \begin{proof} 
Given $\mathbf A, \mathbf  B, \mathbf  C \in \mathcal V_{h,k}$
as in the assumptions \eqref{apeq} of SAP,  first construct
an amalgamating structure 
 $\mathbf  D^r$ for the reducts in the language of $\mathcal V$ as in the proof
of Theorem \ref{thmc}, but here take 
$ \bar{d} \notin A \cup B$, a possibility mentioned in
 Remarks \ref{varr}(b) and \ref{fullc}(a).
Thus here we are taking $D= A \cup B \cup \{  \bar{d} \} $.  

The structure $\mathbf  D^r$ expands in a unique way
to a structure $\mathbf  D$  for the language of $\mathcal V_{h,k}$, under the assumption 
that $\mathbf  D$  extends
$\mathbf A$ and $\mathbf  B$.
The values of $h$ and $k$ on $A \cup B$
are determined by the above request. 
Then set
$ h(\bar{d}) = \bar{d}$ and $ k(\bar{d}) = \bar{d}$.
Such identities should be satisfied if we want $h$ and $k$
to be one the inverse of the other. 
Thus we have a structure $\mathbf  D$ which extends both 
$\mathbf A$ and $\mathbf  B$ and 
$k(h(x))= h(k(x))=x$ are satisfied in $\mathbf  D$,
since they are satisfied in both $\mathbf A$ and $\mathbf  B$. 

It remains to show that
all the requested equations
of the form \eqref{pup} are satisfied in $\mathbf  D$.
To this end, we first observe that, 
for every operation $f$
and every sequence  
$ d_1, d_2, \dots \in A \cup B \cup \bar{d}$,
the value of $f(d_1, d_2, \dots )$
is forced by one of the conditions (E1), (E2)
in the proof of Theorem \ref{thmc}
if and only if  the value of 
$f(e_1, e_2, \dots )$
 is forced by the same condition, where 
$e_1=h(f_1)$,  $e_2=h(f_2)$ \dots \ 
Indeed,  the terms associated to  
$f(d_1, d_2, \dots )$
and  to $f(e_1, e_2, \dots )$ according to (*)
are identical (modulo the ordering of the variables),
since, say, $d_i$ is the interpretation of some constant
if and only if $h(d_i)$  is the interpretation of (the same)
constant, because of the second condition in \eqref{kh}.
Moreover, since $h$ is bijective,
then $d_i = d_j$ if and only if 
$e_i = e_j$.
Thus if the value of 
$f(d_1, d_2, \dots )$ is forced by condition (E1),
then the value of $f(e_1, e_2, \dots )$
is forced by (E1), and conversely. Moreover, if the value forced for 
$f(d_1, d_2, \dots )$ is $d_i$,
then the value forced for  
$f(e_1, e_2, \dots )$ is $e_i=h(d_i)$. 

On the other hand, 
the value of $f(d_1, d_2, \dots )$ 
is forced by (E2)
if and only if 
either $\{ d_{i_1}, d_{i_2}, \dots,
d_{i_{ \ell -1}} \} \subseteq A$
or  $\{ d_{i_1}, d_{i_2}, \dots,
d_{i_{ \ell -1}} \} \subseteq B$,
under the conventions in (E2).
Since the image of $A$ under $h$
is $A$ itself, then  
$\{ d_{i_1}, d_{i_2}, \dots,
d_{i_{ \ell -1}} \} \subseteq A$
if and only if $\{ e_{i_1}, e_{i_2}, \dots,
e_{i_{ \ell -1}} \} =
\{ h(d_{i_1}), h(d_{i_2}), \dots, \allowbreak 
h(d_{i_{ \ell -1}}) \}
 \subseteq A$ and similarly for $B$.
Hence 
the value of $f(d_1, d_2, \dots )$ 
is forced by (E2)
if and only if 
the value of $f(e_1, e_2, \dots )$ 
is forced by (E2).

In the remaining case, 
neither the value of $f(d_1, d_2, \dots )$
nor the value of  $f(e_1, e_2, \dots )$
are forced, hence the proof of \ref{thmc} gives  
$f(d_1, d_2, \dots ) = 
f(e_1, e_2, \allowbreak  \dots )= \bar{d}=h(\bar{d}) $ in $\mathbf  D$.

The proof is almost complete.
We need to show that, limited to those operations $f$
for which \eqref{pup} is assumed in $\mathcal V_{h,k}$,
the very same equation \eqref{pup} holds in $\mathbf  D$.
That is, if  $ d^* =f(d_1, d_2, \dots )$, then
 $ h(d^*) =f(e_1, e_2, \dots )$ in $\mathbf  D$, under the above conventions.
We have already shown in the course of the above arguments
that this identity holds both when the values are forced by condition (E1)
and when the values are not forced.

 It remains to treat the case
when the values are forced by (E2).
In this case the 
values forced for $f(d_1, d_2, \dots )$
and $f(e_1, e_2, \dots )$
are, respectively, 
$t(d_{i_1}, d_{i_2}, \dots, d_{i_{ \ell -1}};d, \dots, d)$
and 
$t(e_{i_1}, e_{i_2}, \dots,
e_{i_{ \ell -1}};
 d, \dots, d)$, evaluated   in $\mathbf A$
or in $\mathbf  B$. For the sake of brevity,
suppose from now on that everything is evaluated in $\mathbf A$.
Let $\hat d_i = d_i$, $\hat e_i = e_i$
if $d_i$ corresponds to an ordinary variable of $t$
and  let $\hat d_i = d$, $\hat e_i = h(d)$
if $d_i$ corresponds to an exceptional variable of $t$.
Then in $\mathbf A$  we have 
\begin{multline}  \labbel{ml}  
t(e_{i_1}, e_{i_2}, \dots, e_{i_{ \ell -1}}; d, \dots, d)=
t(e_{i_1}, e_{i_2}, \dots, e_{i_{ \ell -1}}; h(d), \dots, h(d))=^{def}
\\
f( \hat e_1, \hat e_2, \dots ) =
f( h(\hat d_1), h(\hat d_2), \dots ) =^{\eqref{pup}}
h(f( \hat d_1, \hat d_2, \dots ))=^{def}
\\
h(t(d_{i_1}, d_{i_2}, \dots, d_{i_{ \ell -1}};d, \dots, d)),
\end{multline}   
where the first equality follows 
from the fact that the  variables after the semicolon 
are exceptional
 in the enumeration
of the arguments of $t$, the equalities labeled
with $def$   
follow from the definitions of the term $t$
and of $\hat d_i$, $\hat e_i$ and, finally, we can apply    
\eqref{pup} since  $\hat d_1, \hat d_2, \dots $ 
all belong to  $A$
and $\mathbf A$ satisfies 
\eqref{pup} by assumption.
We have showed that
if equation \eqref{pup} holds in $\mathcal V$
with regard to $f$ and   
$ d^* =f(d_1, d_2, \dots )$, then
 $ h(d^*) =f(e_1, e_2, \dots )= f(h(d_1), h(d_2), \dots )$,
where $ f(d_1, d_2, \dots )$ and
 $ f(e_1, e_2, \dots )$  
are the values forced in $\mathbf  D$ by (E2).
In conclusion equation \eqref{pup}
holds in $\mathbf  D$, hence $\mathbf  D$ 
witnesses SAP. 

If  $\mathcal V$ is axiomatized by a set of equilinear equations
in a language without constants,
then the simpler arguments used in the proof of Theorem \ref{thm} 
apply. In this case we do not need to use elements 
of $C$; compare the comment shortly after 
the introduction of eventuality (E2) in the proof 
of Theorem \ref{thmc}. 
Hence we can repeat the argument 
taking  $C= \emptyset $ and we get JEP. 
As in the proof of Corollary \ref{corthmc},
in the presence of a one-element subalgebra of every algebra
 AP implies JEP.
\end{proof}

As a final remark, the arguments from the present note
are likely to be generalizable to a broader
setting, but we do not know how far.
The main open problem is to ascertain which
arguments can be extended to languages with
further relations besides equality.
In this sense, the main obstacle to generalizations appears to be the
possible mutual incompatibility of AP with conditions corresponding to
\eqref{linec} and \eqref{dc}  
or, put in another way, 
the fact that the outcomes given by conditions similar to (E1) and (E2)
might not agree in all the overlapping cases.

We present a simple example in which
we renounce to the analogue of \eqref{linec}
but we can work with an equilinear analogue of \eqref{dc}.
The example might be a starting point
for possible generalizations.

In the following proposition $\mathcal L$ is an arbitrary 
language, possibly containing relation symbols.

\begin{proposition} \labbel{morel}
Suppose that $T$ is a theory  in the language $\mathcal L$
and $T$ has the (strong) amalgamation
property. Suppose further that every model 
$\mathbf  D_1$  of $T$ can be extended to some model $\mathbf  D$
with an element  $ \bar{d} \in D$
such that $R(\bar{d},\bar{d}, \dots, \bar{d})$
holds in $\mathbf  D$, for every relation symbol 
$R \in \mathcal L$.

Let $T' \supseteq T$ be a theory in a language
   $\mathcal L' \supseteq \mathcal L$ 
such that $\mathcal L' \setminus  \mathcal L$
consists only of function symbols.
Suppose that $T' \setminus  T$
consists only of axioms of the form
\begin{multline}\labbel{relax}
R(f_1(x_{1,1},x_{1,2}, \dots, x_{1,n_1} ), 
f_2(x_{2,1},x_{2,2}, \dots, x_{2,n_2} ), \dots
\\ \dots,
f_m(x_{m,1},x_{m,2}, \dots, x_{m,n_m} )),
   \end{multline}    
for $R \in \mathcal L$, 
 $f_1, f_2, \dots, f_m \in \mathcal L' \setminus  \mathcal L$
and  where
\begin{equation*}    
 \{  x_{1,1},x_{1,2}, \dots, x_{1,n_1}\} =
\{ x_{2,1},x_{2,2}, \dots, x_{2,n_2} \} = \dots =
  \{ x_{m,1},x_{m,2}, \dots, \allowbreak  x_{m,n_m}  \} 
 \end{equation*}
(we are not considering multiplicities.)

Then $T'$ has the (strong) amalgamation property.  
 \end{proposition} 

 \begin{proof} 
Given $\mathbf A'$, $\mathbf  B'$ and $\mathbf  C'$ models of 
$T'$ in  $\mathcal L'$ as in \eqref{apeq},
their reducts to $\mathcal L$  can be amalgamated to 
a model $\mathbf  D$ of $T$, since $T$ satisfies AP.
By the additional assumption on $T$, it is no loss of generality to assume that
there is $ \bar{d} $ in $ D$
such that $R(\bar{d},\bar{d}, \dots, \bar{d})$
holds in $\mathbf  D$ for every relation in $\mathcal L$.

Expand  $\mathbf  D$  to a model $\mathbf  D'$  for $\mathcal L'$ by setting,
for every $f \in \mathcal L' \setminus  \mathcal L$ and 
 $d_1, d_2, \dots \in D$:
 \begin{align*}
 f(d_1, d_2, \dots ) &= f _{\mathbf A} (d_1, d_2, \dots ) 
&&\text{if  $d_1, d_2, \dots \in A$,} 
\\
 f(d_1, d_2, \dots ) &= f _{\mathbf B} (d_1, d_2, \dots )
&&\text{if  $d_1, d_2, \dots \in B$,}
\\
 f(d_1, d_2, \dots ) &= \bar{d}  
&&\text{otherwise}. 
   \end{align*}   

As standard by now, the first two clauses assure that 
$\mathbf  D'$ extends both $\mathbf A'$ and $\mathbf  B'$.
The third  clause and the properties of $ \bar{d}$
imply that each instance of \eqref{relax} is satisfied in $\mathbf  D'$,
provided it is satisfied in $\mathbf A'$, $\mathbf  B'$ and $\mathbf  C'$.
Just notice  that, for any evaluation of \eqref{relax} in $\mathbf  D'$,
the assumption on the occurrences of the variables implies that 
$f_1$ obeys one of the above clauses if and only if 
each $f_i$ obeys the very same clause.  

Hence $\mathbf  D'$ is an amalgamating structure
for $\mathbf A'$, $\mathbf  B'$ and $\mathbf  C'$.  
\end{proof}

\end{document}